\RequirePackage{fix-cm}
\documentclass[11pt,a4paper]{article}

\usepackage[all]{xy}				    
\usepackage{amsthm}				
\usepackage{amsfonts}				
\usepackage{bussproofs}			
\usepackage{cancel}				    
\usepackage{proof}			
\usepackage{geometry}				
\usepackage{graphicx}				
\usepackage{hyperref}				
\usepackage[utf8]{inputenc}		
\usepackage{lineno}				    
\usepackage{mathrsfs}				
\usepackage{MnSymbol}			
\usepackage{hyperref}  
\usepackage{tikz}				 
\usetikzlibrary{tikzmark}
\usepackage{url}				        
\usepackage{xcolor}				    

\usepackage[numbers]{natbib}
  
\theoremstyle{definition}         
\newtheorem{df}{Definition}	   
\newtheorem{cor}[df]{Corollary}   
\newtheorem{lm}[df]{Lemma}	  		
\newtheorem{prop}[df]{Proposition}	 
\newtheorem{rmrk}[df]{Remark}	  
\newtheorem{nt}[df]{Note}		
\newtheorem{tr}[df]{Theorem}	 	

\newcommand{\mc}{\mathcal}

\def\P{\mathbf P}

\usepackage[shortlabels]{enumitem}

\newcommand{\D}{\ensuremath{\mathcal{D}}}

\label{xLOGICS}{}

\newcommand{\lfis}{\textit{LFI}s}

\newcommand{\letf}{$LET_{F}$}

\newcommand{\nletj}{${N}{L}{E}{T}_{J}$}
\newcommand{\gletj}{$GLET_{J}$}
\newcommand{\nletf}{$NLET_{F}$}
\newcommand{\gletf}{$GLET_{F}$}
\newcommand{\letj}{$LET_{J}$}
\newcommand{\fde}{$FDE$}

\newcommand{\nel}{\textit{N4}}
\newcommand{\lets}{\textit{LET}s}

\label{xCOLORS} %

\label{xSYMBOLS}

\newcommand{\cons}{\ensuremath{{\circ}}}
\newcommand{\con}{\ensuremath{{\circ}}}

\newcommand{\To}{\Rightarrow}

\label{xGENERAL COMMANDS}%

\newcommand{\mh}{\noindent}

\newcommand{\m}{\vspace{1mm}}
\newcommand{\mm}{\vspace{2mm}}
\newcommand{\mmm}{\vspace{3mm}}

\newcommand{\bqu}{\begin{quote}} 
\newcommand{\equ}{\end{quote}}
\newcommand{\enr}{\begin{enumerate}[label={(\arabic*)}, resume]}
\newcommand{\eenr}{\end{enumerate}}

\label{xPORTUGUES}

\title{\bf Normalization and cut-elimination theorems for some logics of evidence and truth\thanks{The authors are grateful to N. Kamide for helpful discussions on his paper~\cite{Kam}.  The first author acknowledges support from the National Council for Scientific and Technological Development (CNPq, Brazil), grants  309830/2023-0 and 408040/2021-1. He was also supported by the S\~ao Paulo Research Foundation (FAPESP, Brazil), thematic project {\em Rationality, logic and probability -- RatioLog}, grant  2020/16353-3. 
 The research of the second author was partially supported by the visiting researcher award program funded by 
FAPESP, grant 2022/03862-2.  
 The third author acknowledges support from CNPq, Brazil, grants 310037/2021-2 and 408040/2021-1, and Minas Gerais State Agency for Research and Development (FAPEMIG, Brazil), grant  APQ-02093-21.}}

\author{
    Marcelo E. Coniglio $^\textup{\scriptsize a}$
   \and Martín Figallo $^\textup{\scriptsize b}$ 
   \and Abilio Rodrigues  $^\textup{\scriptsize c}$
}

\date{
    $^\textup{\scriptsize a}$\textit{\small Centre for Logic, Epistemology and the History of Science (CLE), and Institute of Philosophy and the Humanities (IFCH), Universiy of Campinas (UNICAMP). Campinas, Brazil}
    \\
    $^\textup{\scriptsize b}$\textit{\small Departamento de Matem\'atica and Instituto de Matem\'atica (INMABB), Universidad Nacional del Sur (UNS), Bah\'{\i}a Blanca, Argentina.}
    \\
    $^\textup{\scriptsize c}$\textit{\small Department of Philosophy, Federal University of Minas Gerais (UFMG),  Belo Horizonte, Brazil}
}

\begin{document}

\maketitle

\begin{abstract} \mh
In this paper, we investigate proof-theoretic aspects of the logics of evidence and truth
\letj\ and \letf. These logics extend, respectively, Nelson's logic \nel\ and the logic of first-degree 
entailment (\fde), also known as Belnap-Dunn four-valued logic, with a classicality operator \con\ 
that recovers classical 
logic for formulas in its scope. We will present natural deduction and sequent systems for \letj\ 
and \letf, together with 
proofs of normalization and cut-elimination theorems, respectively. As a corollary, we obtain decidability for both logics.

%

\vspace*{4mm}

\end{abstract}

\section{Introduction}

Logics of evidence and truth (\lets) are paracomplete and paraconsistent logics 
that extend the logic of first degree entailment (\fde), also known 
as Belnap-Dunn four-valued logic \cite{belnap.1977.how,dunn76},  
 with a classicality operator \con\ that recovers classical negation for sentences in its scope by means of the following inferences:  
\begin{description} 
\item (1) $ \con A,A,\neg A \vdash  B$, 
\item  (2) $ \con A \vdash A\lor\neg A$. 
\end{description} 
 \mh \lets\  developed from logics of formal inconsistency (\lfis, see \cite{car.mar.2002.tax,car.con.mar.2007}), 
 and both \lets\  and \lfis\ are part of an evolutionary line that begins in the 
 1960s with da Costa's work on paraconsistency  \cite{costa1963}. 
 The intended intuitive interpretation of  \lets\ can be in terms of 
 evidence, which can be either conclusive or non-conclusive, 
or in terms of information. In the latter case, 
the intuitive reading of a formula $\con A$ is that the information conveyed by $A$, positive or negative, is reliable, and 
 when $\con A$ does not hold, it means that there is no reliable information about $A$. 
 \lets, thus, can be interpreted as information-based logics, i.e. logics that are suitable 
 for processing information, in the sense of taking a database as a set of premises 
 and drawing conclusions from those premises in a sensible way, an idea that can be traced back to Belnap-Dunn's papers of the 1970s, in particular to
  \cite[pp.~35-36]{belnap.1977.how}.\footnote{For historical and conceptual 
   aspects of \lets\ see \citep{rod.ant.lu,rod.car.llp}.} 
    
 The logic \letj, introduced in \citep{letj}, extends Nelson's logic \nel\ \citep{almu.nel.84}  
 with the inferences  (1) and (2) above. \nel\ extends \fde\ with a constructive implication: 
 Peirce's law does not hold, but the equivalence between $\neg(A\rightarrow B)$ and $A \land \neg B$ holds. 
 The logic \letf, introduced in \citep{letf},  extends \fde\ with (1) and (2) above plus 
  the  following inferences: 
\begin{description} 
\item (3)  $ \con A,  \neg\con A  \vdash  B$, 
\item  (4) $ B \vdash \con A\lor  \neg\con A$. 
\end{description} 
\mh Natural deduction systems for  \letj\ and \letf\ together with sound and complete 
valuation semantics have been presented in \cite{letj} and \cite{letf}, respectively, 
and Kripke semantics for both logics in \cite{axioms}.

In this paper, we investigate the proof theory of the logics 
\letj\ and \letf. 
Section \ref{sec.nd.letj.letf} presents the natural deduction systems  \nletj\ and \nletf, and in 
Section \ref{sec.normalization} we prove normalization theorems for both systems. 
In Sections \ref{sec.sequent.letj} and \ref{sec.cut.elim} we present a sequent calculus system for \letj\ 
 and the corresponding cut-elimination theorem, and in Section \ref{sec.sequent.letf}
 these results are extended to \letf.

\section{Natural deduction systems} \label{sec.nd.letj.letf}

The language ${\cal L}_J$ of \letj\ is composed of denumerably many sentential letters $p_1, p_2, \dots$, 
the unary connectives $\con$ and $\neg$, the binary connectives $\land$, $\lor$, and $\rightarrow$, 
and parentheses. 
The set of formulas of ${\cal L}_J$, which we will also denote by ${\cal L}_J$, 
is inductively defined in the usual way.  Roman capitals
$A, B,C, \dots$ will be used as metavariables for the formulas of ${\cal L}_J$, 
while Greek capitals $\Gamma, \Delta, \Sigma, \dots$ will be
used as metavariables for sets of formulas. 
 
The proof system \nletj\  will be defined following the notational conventions given in \cite{TS}. Deductions in the system of natural deduction are generated as follows.

\begin{df} (The system \nletj) \ 

\m\mh Deductions in \nletj\ are inductively defined as follows: 

\m\mh Basis: The proof-tree with a single occurrence of an assumption $A$ with a marker is a deduction with conclusion $A$ from open assumption $A$. 

\m\mh Inductive step: Let ${\cal D}$, ${\cal D}_1$, ${\cal D}_2$, ${\cal D}_3$ be deductions. Then, they can be extended by one of the following rules below. The classes {\rm[$\neg A$]$^u$}, {\rm[$\neg B$]$^v$}, {\rm[$A$]$^u$}, {\rm[$B$]$^v$} below contain open assumptions of the deductions of the premises of the final inference, but are closed in the whole deduction. 

{\center

\AxiomC{${\cal D}_1$}
\noLine
\UnaryInfC{$A$}
\AxiomC{${\cal D}_2$}
\noLine
\UnaryInfC{$B$}
\RightLabel{\rm $\wedge$I}
\BinaryInfC{$A \wedge B$}
\DisplayProof \hspace{1cm}
\AxiomC{${\cal D}$}
\noLine
\UnaryInfC{$A\wedge B$}
\RightLabel{\rm $\wedge$E}
\UnaryInfC{$A$}
\DisplayProof 
\AxiomC{${\cal D}$}
\noLine
\UnaryInfC{$A\wedge B$}
\UnaryInfC{$B$}
\DisplayProof 

\mm 

\AxiomC{${\cal D}$}
\noLine
\UnaryInfC{$A$}
\RightLabel{\rm $\vee$I}
\UnaryInfC{$A\vee B$}
\DisplayProof
\AxiomC{${\cal D}$}
\noLine
\UnaryInfC{$B$}
\UnaryInfC{$A\vee B$}
\DisplayProof \hspace{1cm}
\AxiomC{}
\noLine
\UnaryInfC{${\cal D}_1$}
\noLine
\UnaryInfC{$A\vee B$}
\AxiomC{\rm[$A$]$^u$}
\noLine
\UnaryInfC{${\cal D}_2$}
\noLine
\UnaryInfC{$C$}
\AxiomC{\rm [$B$]$^v$}
\noLine
\UnaryInfC{${\cal D}_3$}
\noLine
\UnaryInfC{$C$}
\RightLabel{\rm $\vee$E,$u$,$v$}
\TrinaryInfC{$C$}
\DisplayProof 

\mm 

\AxiomC{\rm [$A$]$^v$}
\noLine
\UnaryInfC{${\cal D}$}
\noLine
\UnaryInfC{$B$}
\RightLabel{\rm $\rightarrow$I,$u$}
\UnaryInfC{$A\rightarrow B$}
\DisplayProof \hspace{1cm}
\AxiomC{${\cal D}_1$}
\noLine
\UnaryInfC{$A\rightarrow B$}
\AxiomC{${\cal D}_2$}
\noLine
\UnaryInfC{$B$}
\RightLabel{\rm $\rightarrow$E}
\BinaryInfC{$B$}
\DisplayProof

\mmm 

\AxiomC{${\cal D}$}
\noLine
\UnaryInfC{$\neg A$}
\RightLabel{\rm $\neg\wedge$I}
\UnaryInfC{$\neg(A\wedge B)$}
\DisplayProof
\AxiomC{${\cal D}$}
\noLine
\UnaryInfC{$\neg B$}
\UnaryInfC{$\neg(A\wedge B)$}
\DisplayProof \hspace{1cm}
\AxiomC{}
\noLine
\UnaryInfC{${\cal D}_1$}
\noLine
\UnaryInfC{$\neg(A\wedge B)$}
\AxiomC{\rm[$\neg A$]$^u$}
\noLine
\UnaryInfC{${\cal D}_2$}
\noLine
\UnaryInfC{$C$}
\AxiomC{\rm [$\neg B$]$^v$}
\noLine
\UnaryInfC{${\cal D}_3$}
\noLine
\UnaryInfC{$C$}
\RightLabel{\rm $\neg\wedge$E,$u$,$v$}
\TrinaryInfC{$C$}
\DisplayProof 

\mm

\AxiomC{${\cal D}_1$}
\noLine
\UnaryInfC{$\neg A$}
\AxiomC{${\cal D}_2$}
\noLine
\UnaryInfC{$\neg B$}
\RightLabel{\rm $\neg\vee$I}
\BinaryInfC{$\neg(A \vee B)$}
\DisplayProof \hspace{1cm}
\AxiomC{${\cal D}$}
\noLine
\UnaryInfC{$\neg(A\vee B)$}
\RightLabel{\rm $\neg\vee$E}
\UnaryInfC{$\neg A$}
\DisplayProof 
\AxiomC{${\cal D}$}
\noLine
\UnaryInfC{$\neg(A\vee B)$}
\UnaryInfC{$\neg B$}
\DisplayProof 

\mm 

\AxiomC{${\cal D}_1$}
\noLine
\UnaryInfC{$A$}
\AxiomC{${\cal D}_2$}
\noLine
\UnaryInfC{$\neg B$}
\RightLabel{\rm $\neg\rightarrow$I}
\BinaryInfC{$\neg(A \rightarrow B)$}
\DisplayProof \hspace{1cm}
\AxiomC{${\cal D}$}
\noLine
\UnaryInfC{$\neg(A\rightarrow B)$}
\RightLabel{\rm $\neg\rightarrow$E}
\UnaryInfC{$A$}
\DisplayProof 
\AxiomC{${\cal D}$}
\noLine
\UnaryInfC{$\neg(A\rightarrow B)$}
\UnaryInfC{$\neg B$}
\DisplayProof

\mm

\AxiomC{${\cal D}$}
\noLine
\UnaryInfC{$A$}
\RightLabel{\rm $\neg\neg$I}
\UnaryInfC{$\neg \neg A$}
\DisplayProof \hspace{1cm}
\AxiomC{${\cal D}$}
\noLine
\UnaryInfC{$\neg\neg A$}
\RightLabel{\rm $\neg\neg$E}
\UnaryInfC{$A$}
\DisplayProof

\mm 

\AxiomC{${\cal D}_1$}
\noLine
\UnaryInfC{$\con A$}
\AxiomC{${\cal D}_3$}
\noLine
\UnaryInfC{$\neg A$}
\AxiomC{${\cal D}_2$}
\noLine
\UnaryInfC{$A$}
\RightLabel{\rm $EXP^\con$}
\TrinaryInfC{B}
\DisplayProof \hspace{1cm}
\AxiomC{}
\noLine
\UnaryInfC{${\cal D}_1$}
\noLine
\UnaryInfC{$\con A$}
\AxiomC{\rm[$A$]$^u$}
\noLine
\UnaryInfC{${\cal D}_2$}
\noLine
\UnaryInfC{$B$}
\AxiomC{\rm [$\neg A$]$^v$}
\noLine
\UnaryInfC{${\cal D}_3$}
\noLine
\UnaryInfC{$B$}
\RightLabel{\rm $PEM^\con$,$u$,$v$}
\TrinaryInfC{$B$}
\DisplayProof 

}
\end{df}

\m\mh 
The notions of {\em major} and {\em minor} premises of the  \con-free rules
are defined as usual. 
$\con A$ is the major premise of $PEM^\con$, the other premises are minor premises. There are no notions of 
 major and minor premises of $EXP^\con$.

\m\mh Let $\Gamma \cup \{A\}$ be a set of formulas of ${\cal L}_J$. We say that the conclusion $A$ is derivable from a set $\Gamma$ of premises,  $\Gamma\vdash_{\mbox{\nletj}} A$, if and only if there is a deduction in \nletj\ with final formula $A$ and whose open assumptions are in $\Gamma$.

\begin{rmrk}\label{rembot}
Note that we can consider the symbol $\bot$ (it behaves here as an arbitrary unprovable propositional constant) as an abbreviation of $\con A\wedge \neg A\wedge A$, for any formula $A$. Then, it is clear that the following rule is derivable in \nletj :
$$ \infer[(\bot)] {B} {\bot}$$
\end{rmrk}

\

Now we turn to the logic \letf. The natural deduction system  \nletf\ can be defined by adding to the implication-free fragment of \letj\ rules corresponding to the inferences (3) and (4) mentioned in the introduction.

\begin{df} (The system \nletf) \ 

\m\mh Let ${\cal L}_F$ be the language obtained from ${\cal L}_J$ by dropping the binary connective 
$\rightarrow$. The calculus \nletf\ results from adding the following rules to the rules of \nletj\  
and dropping the rules $\rightarrow$I and $\rightarrow$E:

\begin{center}
\AxiomC{${\cal D}_1$}
\noLine
\UnaryInfC{$\con A$}
\AxiomC{${\cal D}_2$}
\noLine
\UnaryInfC{$\neg{\con} A$}
\RightLabel{\rm $CONS$}
\BinaryInfC{B}
\DisplayProof \hspace{1.3cm}
\AxiomC{[$\con A$]$^u$}
\noLine
\UnaryInfC{${\mc D}_1$}
\noLine
\UnaryInfC{$B$}
\AxiomC{[$\neg{\con} A$]$^v$}
\noLine
\UnaryInfC{${\mc D}_2$}
\noLine
\UnaryInfC{$B$}
\RightLabel{\rm $COMP, u, v$}
\BinaryInfC{$B$}
\DisplayProof 
\end{center}

\m\mh We will use $\vdash_{\mbox{\letf}}$ to denote the derivability relation generated by \nletf.
\end{df}



\section{Normalization theorems for \nletj\ and \nletf} \label{sec.normalization}


We use the term {\em del-rule} (from ``disjunction-elimination-like'') as in \cite{TS}. We consider rules $\vee$E, $\neg\wedge$E, $PEM^\con$ and $COMP$ as del-rules.

\m

Let $\mc D$ be a deduction in \letj. Recall that a {\em segment} (of length $n$) in the deduction of $\mc D$  is a sequence $A_1, \dots, A_n$ of consecutive occurrences of a formula $A$ in $\mc D$ such that:
\begin{itemize}
\item[(i)] for $1 < n$, $i < n$, $A_i$ is a minor premise of a del-rule application in $\mc D$, with conclusion $A_{i+1}$,
\item[(ii)] $A_n$ is not a minor premise of a del-rule application, 
\item[(iii)] $A_1$ is not the conclusion of a del-rule application. 
\end{itemize}

A segment is {\em maximal}, or a {\em cut (segment)} if $A_n$ is the major premise of an elimination rule (E-rule), and
either $n > 1$, or $n = 1$ and $A_1$ is $A_n$ and $A_n$  is the conclusion of an introduction rule (I-rule) or $EXP^\con$. 
The {\em cutrank} $cr(\sigma)$ of a maximal segment $\sigma$ with formula $A$ is  the complexity of $A$, $|A|$. 
The {\em cutrank} $cr(\mc D)$ of a deduction $\mc D$ is the maximum of the cutranks of cuts in $\mc D$. 
If there is no cut, the cutrank of $\mc D$ is zero. 
A {\em critical cut} of $\mc D$ is a cut of maximal cutrank among all cuts in $\mc D$. We shall use 
$\sigma$, $\sigma'$ to designate segments.
A given deduction $\mc D$ is {\em normal} if it contains no cuts.

\

{\em Detour conversions}

\

We first show how to eliminate cuts of length $1$. The conversions for the connectives $\wedge$, $\vee$ and $\rightarrow$ are the same as those for the respective connectives in the natural deduction system for  intuitionistic logic.

\

{\em $\neg\wedge$-conversion}

\

\begin{center}

\AxiomC{$\mc D_i$}
\noLine
\UnaryInfC{$\neg A_i$}
\UnaryInfC{$\neg(A_1 \wedge A_2)$}

\AxiomC{[$\neg A_1$]}
\noLine
\UnaryInfC{$\mc D'_1$}
\noLine
\UnaryInfC{$C$}

\AxiomC{[$\neg A_2$]}
\noLine
\UnaryInfC{$\mc D'_2$}
\noLine
\UnaryInfC{$C$}

\TrinaryInfC{$C$}
\DisplayProof \hspace{1cm} converts to \hspace{1cm} 
\AxiomC{$\mc D_i$}
\noLine
\UnaryInfC{$\neg A_i$}
\noLine
\UnaryInfC{$\mc D'_i$}
\noLine
\UnaryInfC{$C$}
\DisplayProof
\end{center}

\

{\em $\neg\vee$-conversion}

\

\begin{center}

\AxiomC{$\mc D_1$}
\noLine
\UnaryInfC{$\neg A_1$}

\AxiomC{$\mc D_2$}
\noLine
\UnaryInfC{$\neg A_2$}
\BinaryInfC{$\neg(A_1 \vee A_2)$}
\UnaryInfC{$\neg A_i$}
\DisplayProof \hspace{1cm} converts to \hspace{1cm} 
\AxiomC{$\mc D_i$}
\noLine
\UnaryInfC{$\neg A_i$}
\DisplayProof
\end{center}

\

{\em $\neg\rightarrow$-conversion}

\

\begin{center}

\AxiomC{$\mc D_1$}
\noLine
\UnaryInfC{$A_1$}

\AxiomC{$\mc D_2$}
\noLine
\UnaryInfC{$\neg A_2$}
\BinaryInfC{$\neg(A_1 \rightarrow A_2)$}
\UnaryInfC{$A_1$}
\DisplayProof \hspace{1cm} converts to \hspace{1cm} 
\AxiomC{$\mc D_1$}
\noLine
\UnaryInfC{$A_1$}
\DisplayProof
\end{center}

\ 

and 

\

\begin{center}

\AxiomC{$\mc D_1$}
\noLine
\UnaryInfC{$A_1$}

\AxiomC{$\mc D_2$}
\noLine
\UnaryInfC{$\neg A_2$}
\BinaryInfC{$\neg(A_1 \rightarrow A_2)$}
\UnaryInfC{$\neg A_2$}
\DisplayProof \hspace{1cm} converts to \hspace{1cm} 
\AxiomC{$\mc D_2$}
\noLine
\UnaryInfC{$\neg A_2$}
\DisplayProof
\end{center}

\

{\em $EXP^\con$-conversion} 

\

\begin{center}
\AxiomC{$\mc D_1$}
\noLine
\UnaryInfC{$\con A$}
\AxiomC{$\mc D_2$}
\noLine
\UnaryInfC{$A$}
\AxiomC{$\mc D_3$}
\noLine
\UnaryInfC{$\neg A$}
\TrinaryInfC{$B$}

\AxiomC{$\mc D'$}
\RightLabel{\small (E-rule)}
\BinaryInfC{$C$}

\DisplayProof \hspace{0cm} converts to \hspace{0cm} 
\AxiomC{$\mc D_1$}
\noLine
\UnaryInfC{$\con A$}
\AxiomC{$\mc D_2$}
\noLine
\UnaryInfC{$A$}
\AxiomC{$\mc D_3$}
\noLine
\UnaryInfC{$\neg A$}
\RightLabel{\small (E-rule)}
\TrinaryInfC{$C$}
\DisplayProof
\end{center}

\

{\em $CONS$-conversion} 

\

\begin{center}
\AxiomC{$\mc D_1$}
\noLine
\UnaryInfC{$\con A$}
\AxiomC{$\mc D_2$}
\noLine
\UnaryInfC{$\neg{\con} A$}
\RightLabel{\small (CONS)}
\BinaryInfC{$B$}

\AxiomC{$\mc D'$}
\RightLabel{\small (E-rule)}
\BinaryInfC{$C$}

\DisplayProof \hspace{0cm} converts to \hspace{0cm} 
\AxiomC{$\mc D_1$}
\noLine
\UnaryInfC{$\con A$}
\AxiomC{$\mc D_2$}
\noLine
\UnaryInfC{$\neg{\con} A$}
\RightLabel{\small (CONS)}
\BinaryInfC{$C$}
\DisplayProof
\end{center}

\

{\em Permutation conversions}

\

In order to remove cuts of length $> 1$, we permute $E$-rules upwards over minor premises of $\vee$E and $\neg\wedge$E. The $\vee$-permutation conversion is the same as the respective permutation conversion in the intuitionistic logic. 

\

{\em $\vee$-permutation conversion}

\

{\scriptsize
\AxiomC{$\mc D$}
\noLine
\UnaryInfC{$A_1 \vee A_2$}

\AxiomC{[$A_1$]}
\noLine
\UnaryInfC{$\mc D_1$}
\noLine
\UnaryInfC{$C$}

\AxiomC{[$ A_2$]}
\noLine
\UnaryInfC{$\mc D_2$}
\noLine
\UnaryInfC{$C$}

\TrinaryInfC{$C$}

\AxiomC{$\mc D'$}
\RightLabel{\scriptsize (R)}
\BinaryInfC{$B$}
\DisplayProof c. to 
\AxiomC{$\mc D$}
\noLine
\UnaryInfC{$A_1 \vee A_2$}
\AxiomC{[$A_1$]}
\noLine
\UnaryInfC{$\mc D_1$}
\noLine
\UnaryInfC{$C$} 
\AxiomC{$\mc D'$}
\RightLabel{\scriptsize (R)}
\BinaryInfC{$B$}
\AxiomC{[$A_2$]}
\noLine
\UnaryInfC{$\mc D_2$}
\noLine
\UnaryInfC{$C$} \AxiomC{$\mc D'$}
\RightLabel{\scriptsize (R)}
\BinaryInfC{$B$}
\TrinaryInfC{$B$}
\DisplayProof
}

\

{\em $\neg \wedge$-permutation conversion}

\

{\tiny
\AxiomC{$\mc D$}
\noLine
\UnaryInfC{$\neg(A_1 \wedge A_2)$}

\AxiomC{[$\neg A_1$]}
\noLine
\UnaryInfC{$\mc D_1$}
\noLine
\UnaryInfC{$C$}

\AxiomC{[$\neg A_2$]}
\noLine
\UnaryInfC{$\mc D_2$}
\noLine
\UnaryInfC{$C$}

\TrinaryInfC{$C$}

\AxiomC{$\mc D'$}
\RightLabel{\scriptsize (R)}
\BinaryInfC{$B$}
\DisplayProof c. to 
\AxiomC{$\mc D$}
\noLine
\UnaryInfC{$\neg(A_1 \wedge A_2)$}
\AxiomC{[$\neg A_1$]}
\noLine
\UnaryInfC{$\mc D_1$}
\noLine
\UnaryInfC{$C$} 
\AxiomC{$\mc D'$}
\RightLabel{\scriptsize (R)}
\BinaryInfC{$B$}
\AxiomC{[$\neg A_2$]}
\noLine
\UnaryInfC{$\mc D_2$}
\noLine
\UnaryInfC{$C$} \AxiomC{$\mc D'$}
\RightLabel{\scriptsize (R)}
\BinaryInfC{$B$}
\TrinaryInfC{$B$}
\DisplayProof
}

\

{\em $PEM^\con$-conversion} 

\

{\small
\AxiomC{$\mc D$}
\noLine
\UnaryInfC{$\con A$}

\AxiomC{[$\neg A$]}
\noLine
\UnaryInfC{$\mc D_1$}
\noLine
\UnaryInfC{$C$}

\AxiomC{[$ A$]}
\noLine
\UnaryInfC{$\mc D_2$}
\noLine
\UnaryInfC{$C$}

\TrinaryInfC{$C$}

\AxiomC{$\mc D'$}
\RightLabel{\small (R)}
\BinaryInfC{$B$}
\DisplayProof c. to 
\AxiomC{$\mc D$}
\noLine
\UnaryInfC{$\con A$}
\AxiomC{[$\neg A$]}
\noLine
\UnaryInfC{$\mc D_1$}
\noLine
\UnaryInfC{$C$} 
\AxiomC{$\mc D'$}
\RightLabel{\small (R)}
\BinaryInfC{$B$}
\AxiomC{[$A$]}
\noLine
\UnaryInfC{$\mc D_2$}
\noLine
\UnaryInfC{$C$} \AxiomC{$\mc D'$}
\RightLabel{\small (R)}
\BinaryInfC{$B$}
\TrinaryInfC{$B$}
\DisplayProof
}

\

{\em $COMP$-conversion} 

\

{\small

\AxiomC{[$\con A$]}
\noLine
\UnaryInfC{$\mc D_1$}
\noLine
\UnaryInfC{$C$}

\AxiomC{[$\neg{\con} A$]}
\noLine
\UnaryInfC{$\mc D_2$}
\noLine
\UnaryInfC{$C$}

\BinaryInfC{$C$}

\AxiomC{$\mc D'$}
\RightLabel{\small (R)}
\BinaryInfC{$B$}
\DisplayProof converts to 
\AxiomC{[$\con A$]}
\noLine
\UnaryInfC{$\mc D_1$}
\noLine
\UnaryInfC{$C$}
\AxiomC{$\mc D'$}
\RightLabel{\small (R)}
\BinaryInfC{$B$}
\AxiomC{[$\neg{\con} A$]}
\noLine
\UnaryInfC{$\mc D_1$}
\noLine
\UnaryInfC{$C$}
\AxiomC{$\mc D'$}
\RightLabel{\small (R)}
\BinaryInfC{$B$}
\BinaryInfC{$B$}
\DisplayProof
}

\

{\em Simplification conversions}

\

As usual, applications of $\vee$E with major premise $A_1\vee A_2$, where at least one of $[A_1]$,$ [A_2]$
is empty in the deduction of the first or second minor premise, are redundant. Similarly, an application
of $\neg\wedge$E, $PEM^\con$ with major premise  $\neg(A_1\wedge A_2)$, $\con A$ (respectively) where at least one of $[\neg A_1]$ ,$ [\neg A_2]$ and  $[\neg A]$, $[A]$, (respectively) is empty in the deduction of the first or second minor premise, is redundant. Redundant applications
of $\vee$E can be removed as in the intuitionistic case. Similarly, redundant applications of $\neg\wedge$E can be remove easily as follows

\

\begin{center}
\AxiomC{$\mc D$}
\noLine
\UnaryInfC{$\neg(A_1 \wedge A_2)$}

\AxiomC{[$\neg A_1$]}
\noLine
\UnaryInfC{$\mc D_1$}
\noLine
\UnaryInfC{$C$}

\AxiomC{[$\neg A_2$]}
\noLine
\UnaryInfC{$\mc D_2$}
\noLine
\UnaryInfC{$C$}

\TrinaryInfC{$C$}
\DisplayProof \hspace{1cm} converts to \hspace{1cm} 
\AxiomC{$\mc D_i$}
\noLine
\UnaryInfC{$C$}
\DisplayProof
\end{center}
where no assumptions are discharged by $\neg\wedge$E in $\mc D_i$. The other cases are treated similarly.


\

\begin{lm}  \ \label{EXP} 

\m\mh Let $\mc D$ be a deduction of $A$ from $\Gamma$ (in \letj). Then, there is a deduction $\mc D'$ of $A$ from $\Gamma$ in which the consequence of every application of $EXP^\con$ is either a literal (i.e. an atomic formula or its negation) or a formula of the form $\con B$ or $\neg{\con} B$. 
\end{lm}
\begin{proof} Let $c$ be the highest complexity of a consequence of an application of $EXP^\con$ in $\mc D$.  Let $B$ be the consequence of an application of $EXP^\con$ of complexity $c$ and that no consequence of an application of $EXP^\con$ en $\mc D$ that stands above $B$, is  of complexity $c$. Then, $\mc D$ has the form

\begin{prooftree}
\AxiomC{$\mc D_1$}
\noLine
\UnaryInfC{$\con A$}

\AxiomC{$\mc D_2$}
\noLine
\UnaryInfC{$A$}

\AxiomC{$\mc D_3$}
\noLine
\UnaryInfC{$\neg A$}
\TrinaryInfC{$B$}

\AxiomC{$\mc D'$}
\BinaryInfC{$\mc D''$}
\end{prooftree}

Suppose that $B$ is not a literal and that it is not of the form $\con C$ nor $\neg\con C$. Then, $B$ is of the form $C\wedge D$,  $C\vee D$, $C\rightarrow D$, $\neg(C\wedge D)$, $\neg(C\vee D)$, $\neg \neg C$. We remove this application of $EXP^\con$ easily as follows. 

If $B$ is $C\wedge D$ the

\AxiomC{$\mc D_1$}
\noLine
\UnaryInfC{$\con A$}

\AxiomC{$\mc D_2$}
\noLine
\UnaryInfC{$A$}

\AxiomC{$\mc D_3$}
\noLine
\UnaryInfC{$\neg A$}
\TrinaryInfC{$C\wedge D$}
\noLine
\UnaryInfC{$\mc D_4$}
\DisplayProof \hspace{1cm} converts to \hspace{1cm}  \AxiomC{$\mc D_1$}
\noLine
\UnaryInfC{$\con A$}
\AxiomC{$\mc D_2$}
\noLine
\UnaryInfC{$A$}
\AxiomC{$\mc D_3$}
\noLine
\UnaryInfC{$\neg A$}
\TrinaryInfC{$C$}
\AxiomC{$\mc D_1$}
\noLine
\UnaryInfC{$\con A$}
\AxiomC{$\mc D_2$}
\noLine
\UnaryInfC{$A$}
\AxiomC{$\mc D_3$}
\noLine
\UnaryInfC{$\neg A$}
\TrinaryInfC{$D$}
\RightLabel{$\wedge$I}
\BinaryInfC{$C\wedge D$} 
\noLine
\UnaryInfC{$\mc D_4$}
\DisplayProof

\

\

\AxiomC{$\mc D_1$}
\noLine
\UnaryInfC{$\con A$}
\AxiomC{$\mc D_2$}
\noLine
\UnaryInfC{$A$}
\AxiomC{$\mc D_3$}
\noLine
\UnaryInfC{$\neg A$}
\TrinaryInfC{$\neg(C\wedge D)$}
\noLine
\UnaryInfC{$\mc D_4$}
\DisplayProof \hspace{1cm} converts to \hspace{1cm} \AxiomC{$\mc D_1$}
\noLine
\UnaryInfC{$\con A$}
\AxiomC{$\mc D_2$}
\noLine
\UnaryInfC{$A$}
\AxiomC{$\mc D_3$}
\noLine
\UnaryInfC{$\neg A$}
\TrinaryInfC{$\neg C$}
\RightLabel{$\neg\wedge$I}
\UnaryInfC{$\neg(C\wedge D)$}
\noLine
\UnaryInfC{$\mc D_4$}
\DisplayProof 

\

\

\AxiomC{$\mc D_1$}
\noLine
\UnaryInfC{$\con A$}
\AxiomC{$\mc D_2$}
\noLine
\UnaryInfC{$A$}
\AxiomC{$\mc D_3$}
\noLine
\UnaryInfC{$\neg A$}
\TrinaryInfC{$C\rightarrow D$}
\noLine
\UnaryInfC{$\mc D_4$}
\DisplayProof \hspace{1cm} converts to \hspace{1cm} \AxiomC{$\mc D_1$}
\noLine
\UnaryInfC{$\con A$}
\AxiomC{$\mc D_2$}
\noLine
\UnaryInfC{$A$}
\AxiomC{$\mc D_3$}
\noLine
\UnaryInfC{$\neg A$}
\TrinaryInfC{$D$}
\RightLabel{$\rightarrow$I}
\UnaryInfC{$C\rightarrow D$}
\noLine
\UnaryInfC{$\mc D_4$}
\DisplayProof 

\

\mh In this last proof the set of assumptions [$C$] is empty. The rest of the cases are similar. So, the new applications of $EXP^\con$ that arise from these transformations have consequence with complexity $< c$. Thus, by successively repeating the transformation we 
obtain a proof fulfilling the desired properties.
\end{proof}


\begin{tr}\label{trLetj} (Normalization for \letj) 

\m\mh Each derivation $\mc D$ in \letj\ reduces to a normal derivation.
\end{tr}
\begin{proof} Let $\mc D$ be a deduction fulfilling the conditions of Lemma \ref{EXP} and such that every redundant application of any del-rule as well as detours with maximal formula  has been removed as indicated in the {\em simplification conversions} section.   As usual, the proof is by double induction on the cutrank $n$ of $\mc D$ and the sum of lengths of all critical cuts, $m$,  in $\mc D$. By a suitable choice of the critical cut to which we apply a conversion we can achieve that either $n$ decreases, or that $n$ remains constant but $m$ decreases. 

We say that $\sigma$ is a {\em top critical cut} in $\mc D$ if no critical cut occurs in a branch of $\mc D$, above $\sigma$.
Let $\sigma$ be the rightmost top critical cut of $\mc D$ with formula $A$. 
If the length of $\sigma$ is $1$ and it is the only maximal segment in $\mc D$, we apply a conversion to $\sigma$ in $\mc D$ obtaining a new deduction $\mc D'$ that has a lower cutrank. If the length of $\sigma$ is $1$ but it is not the only maximal segment of $\mc D$, then the cutrank of $\mc D'$ is equal to the cutrank of $\mc D$ but the sum of lengths of all critical cuts in $\mc D'$ is lower than the same sum in $\mc D$.
If the length of $\sigma$ is $>1$, we apply a permutation conversion to $\mc D$ obtaining the deduction $\mc D'$ which has the same cutrank as $\mc D$, but with a lower value for $m$.
\end{proof}

\

Similarly, it is possible to eliminate cuts of any derivation in \letf. 

\begin{tr} (Normalization for \nletf) \ 

\m\mh Each derivation $\mc D$ in \letf\ reduces to a normal derivation.
\end{tr}
\begin{proof} We proceed as in the proof of Theorem \ref{trLetj}. Details are left to the reader. 
\end{proof}

\section{Sequent calculus for \letj} \label{sec.sequent.letj}

 Let us start by considering the following single-conclusion sequent calculus, which is a sort of `direct translation' of the rules of \nletj\ into sequent rules,  where the rule $EXP^\con$ of \nletj\ corresponds to an axiom with the same name. 
As usual $\Gamma, \Delta, $ etc. are finite sets of formulas. Besides, we shall identify the sequent \, $\Gamma\To $ \, with the sequent \, $\Gamma\To \bot$.

\begin{df} \ 

\m\mh The sequent calculus $GB$ consists of the following axiom and rules:

\mm

\textit{Axioms:}
\begin{center}
$\infer[Id] {A\To A} {}$
\quad \quad 
${\infer[EXP^\con] {\cons A, A, \neg A, \Gamma\To C} {}} $
\end{center}

\textit{Structural rules:}
\begin{center}
$\infer[LW] {A, \Gamma \To C} {\Gamma \To C}$
\quad 
$\infer[RW] {\Gamma \To A} { \Gamma\To}$

\

$\infer[Cut] {\Gamma \To C} { \Gamma \To A & A, \Gamma \To C}$
\end{center}

\textit{Logical rules:}
\begin{center}
$\infer[L\land] {A \land B,\Gamma \To C} {A,B , \Gamma \To C}$  
\quad 
$\infer[R\land] {\Gamma \To A \land B} { \Gamma\To A   &   \Gamma\To B}$

\

$\infer[L\lor] {A \lor B,\Gamma \To C} {A,\Gamma \To C & B,\Gamma \To C}$
\quad 
$\infer[R\lor] {\Gamma \To A \lor B } { \Gamma\To A }$  $ \infer[] {\Gamma \To A \lor B } { \Gamma\To B}$

\

$\infer[L\neg\lor] {\neg(A \lor B),\Gamma \To C} {\neg A, \neg B,\Gamma \To C}$
\quad 
$\infer[R\neg\lor] {\Gamma \To \neg(A \lor B)} { \Gamma\To \neg A   &   \Gamma\To \neg B}$

\

$\infer[L\to] {A \to B,\Gamma \To C} {\Gamma \To A & B,\Gamma \To C}$
\quad 
$\infer[R\to] {\Gamma \To A \to B} { \Gamma,A \To B}$

\

$\infer[L\neg\to] {\Gamma,\neg(A \to B) \To C} { \Gamma,A,\neg B\To C}$
\quad 
$\infer[R\neg\to] {\Gamma \To \neg(A \to B)} {\Gamma \To A & \Gamma \To \neg B}$

\

$\infer[L\neg\land] {\neg(A \land B),\Gamma \To C} 
{\neg A,\Gamma \To C & \neg B,\Gamma \To C}$
\quad 
$\infer[R\neg\land] {\Gamma \To \neg(A \land B)} {\Gamma\To \neg A}$ 
$\infer[] {\Gamma \To \neg(A \land B)} { \Gamma\To \neg B}$

\

$\infer[PEM^\con] {\cons A, \Gamma \To C}{A, \Gamma\To C & \neg A, \Gamma\To C }$

\

$\infer[L\neg\neg] {\neg\neg A , \Gamma \To B}{A , \Gamma \To B}$
\quad 
$\infer[R\neg\neg] {\Gamma \To \neg\neg A} {\Gamma \To A}$

\end{center}
\end{df}

\mh  We will see that  the system $GB$, although equivalent to \nletj, does not enjoy the cut-elimination property.

\begin{nt} If $G$ is a sequent calculus and $(r)$ is an arbitrary sequent rule, we shall denote by $G \bigplus (r)$ the sequent calculus obtained by adding the rule $(r)$ to the rules of $G$. Similarly, we denote $G \bigplus (r_1)\bigplus \dots \bigplus (r_n)$ the calculus obtained from $G$ by adding the rules $(r_1), \dots (r_n)$. Similarly, if $(r)$ is a rule (or axiom) of $G$ we denote $G - (r)$ the system obtained from $G$ by dropping $(r)$.
\end{nt}


\noindent Then, 


\begin{lm}\label{comp}  \

\m\mh Let $\Gamma \cup \{C\}$ a set of formula. The following conditions are equivalent:
\begin{itemize} 
\item[(i)] $\Gamma \To C$ is provable in $GB$,
\item[(ii)] there is a deduction of $C$ from $\Gamma$ in \nletj .
\end{itemize}
\end{lm}
\begin{proof} (i) implies (ii). Let $\P$ a proof of $\Gamma \To\Delta$ in $GB$. We shall use induction on the number $n$ of rule applications in $\P$, $n\geq 0$.\\
If $n=0$ then we have that $\Gamma\To\Delta$ is (1) $A\To A$ or (2) $\con A, \neg A, A, \Gamma\To C$. In case (1), it is clear that $A \vdash_{\mbox{\nletj} } A$. In case (2),

\begin{prooftree}
\AxiomC{$\con A$}
\AxiomC{$\neg A$}
\AxiomC{$A$}
\RightLabel{$EXP^\con$}
\TrinaryInfC{$C$}
\end{prooftree}

\noindent is a deduction of $C$ from $\{\con A, \neg A, A\}$. Then, using the rules $\wedge$I and $\wedge$E, we can construct a deduction of $C$ from $\Gamma \cup \{\con A, \neg A, A\}$.\\[2mm]
Now, (I.H.) suppose that ``(i)$\To$(ii)'' holds for $n<k$, $k\geq 0$. Let $n=k$, that is $\P$ is a derivation in \gletj\  with last rule $(r)$ of the form
\begin{prooftree}
\AxiomC{$\P$}
\noLine
\UnaryInfC{$\vdots$}
\RightLabel{\small ($r$)}
\UnaryInfC{$\Gamma\Rightarrow\Delta$}
\end{prooftree}
\noindent If $(r)$ is left weakening, then the last rule of $\P$ has the form $\displaystyle \frac{\Gamma'\Rightarrow C}{\Gamma', A\Rightarrow C} (r)$. By (I.H.), there exists a deduction $\cal D$ of $C$ from $\Gamma'$, then 
\begin{prooftree}
\AxiomC{$\cal D$}
\noLine
\UnaryInfC{$C$}
\AxiomC{$A$}
\RightLabel{\small $\wedge$I}
\BinaryInfC{$C \wedge A$}
\RightLabel{\small $\wedge$E}
\UnaryInfC{$C$}
\end{prooftree}
\noindent is a deduction of $C$ from $\Gamma'\cup\{A\}$. If $(r)$ is right weakening, then $(r)$ has the form $\displaystyle \frac{\Gamma\Rightarrow}{\Gamma\Rightarrow C}  (r)$, then by (I.H.) there is a deduction $\cal D$ of $\bot$ from $\Gamma$. Then, using $\wedge$E and $EXP^\con$, we can construct the desired deduction.

If $(r)$ is the cut rule, then $(r)$ is \, $\displaystyle \frac{\Gamma\To A \hspace{.5cm}  A, \Gamma\To B } {\Gamma \To B} \, {cut}$. By (I.H.), there are deductions ${\cal D}_1$ and ${\cal D}_2$ of $A$ from $\Gamma$ and $B$ from $\Gamma \cup \{A\}$, respectively. Then 

\begin{prooftree}
\AxiomC{${\cal D}_1$}
\noLine
\UnaryInfC{$A$}
\noLine
\UnaryInfC{${\cal D}_2$}
\noLine
\UnaryInfC{$B$}
\end{prooftree}
is a deduction of $B$ from $\Gamma$. \\[1.5mm]
If $(r)$ is any of the logical rules $R\vee$, $L\vee$, $R\wedge$, $L\wedge$, $R\to$, $L\to$, $R\vee$, $R\neg\vee$, $L\neg\vee$, $R\neg\wedge$, $L\neg\wedge$, $R\neg\to$ and $L\neg\to$; it is not difficult to find the desired deduction. We shall show it just for $L\neg\to$ and $R\neg\to$. If the last rule application is $\displaystyle \frac{\Gamma, A, \neg B \To C}{\Gamma, \neg(A\to \neg B) \To C}\, {\small L\neg\to}$. Then, by (I.H), we have a deduction $\cal D$ of $C$ from $\Gamma\cup\{A,\neg B\}$. Then

\begin{prooftree}
\AxiomC{$\neg(A\to B)$}
\RightLabel{$\neg\to$E}
\UnaryInfC{$A$}
\AxiomC{$\neg(A\to B)$}
\RightLabel{$\neg\to$E}
\UnaryInfC{$\neg B$}
\BinaryInfC{${\cal D}$}
\noLine
\UnaryInfC{$C$}
\end{prooftree}

\noindent is a deduction of $C$ from $\Gamma\cup\{\neg(A\to \neg B)\}$.\\
Finally, suppose that $(r)$ is $PEM^\con$, $\displaystyle \frac{A, \Gamma\To C \hspace{.5cm} \neg A, \Gamma\To C } {\cons A, \Gamma \To C}\, PEM^\con$ . Then, by (I.H.), we have deductions ${\cal D}_1$ and  ${\cal D}_2$ of $C$ from $\Gamma\cup\{A\}$ and $C$ from $\Gamma\cup\{\neg A\}$. Then

\begin{prooftree}
\AxiomC{$\con A$}

\AxiomC{$[A]^u$}
\noLine
\UnaryInfC{${\cal D}_1$}
\noLine
\UnaryInfC{$C$}

\AxiomC{$[\neg A]^v$}
\noLine
\UnaryInfC{${\cal D}_2$}
\noLine
\UnaryInfC{$C$}

\RightLabel{\small $PEM^\con$, $u$, $v$}
\TrinaryInfC{$C$}
\end{prooftree}

\noindent is a deduction of $C$ from $\Gamma\cup\{\con A\}$. \\

\noindent (ii) implies (i).   Let $\cal D$ be a deduction of $C$ from $\Gamma$ in \nletj .  As before, we use induction on the number $n$ of rule instances in the deduction $\cal D$. If $n=0$ the proof is trivial. (I.H.) Suppose that ``(ii) implies (i)'' holds for $n<k$, $k>0$; and let $(r)$ the last rule instance in $\cal D$. Suppose that $(r)$ is one of the introduction/elimination rule of $\wedge$I, $\wedge$E, $\neg\wedge$I, $\neg\wedge$E, $\vee$I, $\vee$E , $\neg\vee$I, $\neg\vee$E, $\to$I, $\to$E, $\neg\to$I, $\neg\to$E, $\neg\neg$I or $\neg\neg$E. We shall show the construction of the desired sequent proof just for $\neg\wedge$E

$$\infer {C}{  \deduce {\neg (A \land B)} {\D_3} & { \deduce{C}{\deduce {\D_1}{[\neg A]}} } & {\deduce{C}{\deduce {\D_2}{[\neg B]}} } }$$

\noindent By (I.H.), the sequents $\Gamma_1\To \neg(A\wedge B)$, $\Gamma_2, \neg A\To C$ and $\Gamma_3, \neg B\To C$ are probable in $GB$ where $\Gamma_1\cup\Gamma_2\cup\Gamma_3 =\Gamma$. Then, $\Gamma \To C$ is also probable. Indeed, 

$$\infer[cut] {\Gamma \To C}{\infer[LW] {\Gamma \To \neg (A \land B)} {\Gamma_1 \To \neg (A \land B)} & \infer[L\neg\wedge] {  { \Gamma, \neg (A \land B) \To C  }       }{   \infer[LW] { \Gamma, \neg A \To C }  {\Gamma_2, \neg A \To C }  &  \infer[LW]  { \Gamma, \neg B \To C }  {\Gamma_3, \neg B \To C }  }}$$

\noindent  If $(r)$ is $EXP^\con$, then $\cal D$ is of the form 

\begin{prooftree}
\AxiomC{${\cal D}_1$}
\noLine
\UnaryInfC{$\con A$}
\AxiomC{${\cal D}_2$}
\noLine
\UnaryInfC{$A$}
\AxiomC{${\cal D}_3$}
\noLine
\UnaryInfC{$\neg A$}
\RightLabel{\rm $EXP^\con$}
\TrinaryInfC{C}
\end{prooftree}

\noindent By (I.H.), we know that the sequents $\Gamma_1\To\con A$, $\Gamma_2\To \neg A$ and $\Gamma_3\To A$ are probable in $GB$, where $\Gamma_1\cup\Gamma_2\cup\Gamma_3=\Gamma$. Then, $\Gamma\To C$ is also probable in $GB$. Indeed, 

\begin{prooftree}
\AxiomC{$\Gamma_3\To A$}
\UnaryInfC{$\Gamma\To A, C$}

\AxiomC{$\Gamma_2\To \neg A$}
\UnaryInfC{$ A, \Gamma\To \neg A, C$}

\AxiomC{$\Gamma_1\To \con A$}
\UnaryInfC{$\neg A, A, \Gamma\To \con A, C$}

\AxiomC{}
\RightLabel{ $EXP^\con$}
\UnaryInfC{$\con A, \neg A, A, \Gamma \To C$}

\RightLabel{ $cut$}
\BinaryInfC{$\neg A, A, \Gamma \To C$}
\RightLabel{ $cut$}
\BinaryInfC{$A, \Gamma \To C$}
\RightLabel{ $cut$}
\BinaryInfC{$\Gamma \To C$}
\end{prooftree}
If $(r)$ is $PEM^\con$, then the proof goes similarly to the previous case.
\end{proof}

\begin{rmrk}\label{rem1} \ 

\m\mh The system $GB$ does not enjoy the cut-elimination property. Indeed, the following probable sequent does not have a cut-free proof:
$$\con(p\land q), p, q, \neg p\Rightarrow$$

\noindent  where $p$ and $q$ are propositional variables. 
\begin{prooftree}\scriptsize
\AxiomC{}
\RightLabel{$Id$}
\UnaryInfC{$p\To p$}
\RightLabel{$LW$}
\UnaryInfC{$p, q\To p$}

\AxiomC{}
\RightLabel{$Id$}
\UnaryInfC{$q\To q$}
\RightLabel{$LW$}
\UnaryInfC{$p, q\To q$}
\RightLabel{$R\land$}
\BinaryInfC{$p,q\To p\land q$}
\RightLabel{$LW$}
\UnaryInfC{$\con(p\land q), p, q, \neg p\To  p\land q$}

\AxiomC{}
\RightLabel{$Id$}
\UnaryInfC{$\neg p\To \neg p$}
\RightLabel{$R\neg\land$}
\UnaryInfC{$\neg p\To \neg (p\land q)$}
\RightLabel{$LW's$}
\UnaryInfC{$\con(p\land q), p\land q,\neg p\To \neg (p\land q)$}

\AxiomC{}
\RightLabel{$EXP^\con$}
\UnaryInfC{$\con(p\land q), p\land q,\neg(p\land q) \To$}
\RightLabel{$RW$}
\UnaryInfC{$\con(p\land q), p\land q,\neg p, \neg(p\land q) \To$}
\RightLabel{$cut$}
\BinaryInfC{$\con(p\land q), p\land q,\neg p\To$}
\RightLabel{$LW$}
\UnaryInfC{$\con(p\land q), p\land q, p, q, \neg p\To$}
\RightLabel{$LW$}
\RightLabel{$cut$}
\BinaryInfC{$\con(p\land q), p, q,\neg p\To$}
\end{prooftree}
\noindent The last rule in any proof of it has to be an instance of cut.
\end{rmrk}

 Now, we replace the sequent rule $EXP^\con$ with the rule $EXP_1^\con$ below, 
 obtaining a system equivalent to $GB$, 
which is equivalent to \nletj, but enjoys cut-elimination. 

\begin{prop}\label{propeq} \ 

\m\mh Let $G$ be a sequent calculus containing the structural rules of (left and right) weakening and the cut rule and consider the rule
\begin{center}
$\infer[EXP_1^\con] {\con A, \Gamma\To}{\Gamma\To A & \Gamma\To \neg A}$
\end{center}
Then,  $EXP^\con$ is derivable in $G \bigplus EXP_1^\con$ \, and \,  $EXP_1^\con$ is derivable in $G \bigplus EXP^\con$.
\end{prop}
\begin{proof}

\

\begin{prooftree}
\AxiomC{}
\RightLabel{$Id$}
\UnaryInfC{$A \To A$}
\RightLabel{$LW$}
\UnaryInfC{$\neg A, A \To A$}
\AxiomC{}
\RightLabel{$Id$}
\UnaryInfC{$\neg A \To \neg A$}
\RightLabel{$LW$}
\UnaryInfC{$\neg A, A \To \neg A$}
\RightLabel{$EXP_1^\con$}
\BinaryInfC{$\con A, \neg A, A \To $}
\end{prooftree}

\

\begin{prooftree}
\AxiomC{$\Gamma \To \neg A$}
\RightLabel{$LW$}
\UnaryInfC{$\con A, \Gamma \To \neg A$}
\AxiomC{$\Gamma \To A$}
\RightLabel{$LW$}
\UnaryInfC{$\con A, \neg A, \Gamma \To A$}
\AxiomC{}
\RightLabel{$EXP^\con$}
\UnaryInfC{$\con A, \neg A, A,  \Gamma \To $}
\RightLabel{$cut, A$}
\BinaryInfC{$\con A, \neg A,  \Gamma \To $}
\RightLabel{$cut, \neg A$}
\BinaryInfC{$\con A, \Gamma \To $}
\end{prooftree}
\end{proof}

\

\begin{df} (Sequent calculus for \letj)  

\m\mh The sequent calculus \gletj\ is the calculus $\big(GB - EXP^\con  \big) \bigplus EXP_1^\con$. That is, \gletj\ is the system obtained from $GB$ by dropping the axiom $EXP^\con$ and adding the rule:
\begin{center}
$\infer[EXP_1^\con] {\con A, \Gamma\To}{\Gamma\To A & \Gamma\To \neg A}$
\end{center}
\end{df}

\noindent Therefore, we have the following:

\begin{cor} 

\m\mh Let $\Gamma \cup \{C\}$ a set of formula. The following conditions are equivalent:
\begin{itemize}
\item[(i)] the sequent $\Gamma \To C$ is probable in \gletj,
\item[(ii)] there is a deduction of the $C$ from $\Gamma$ in \nletj .
\end{itemize}
\end{cor}

\

\begin{rmrk} The sequent considered in Remark \ref{rem1} has a cut-free proof  for \break $\con(p\land q), p, q, \neg p\Rightarrow$ in \gletj . Indeed,  

\begin{prooftree}\small
\AxiomC{}
\RightLabel{$Id$}
\UnaryInfC{$p\To p$}
\RightLabel{$LW$}
\UnaryInfC{$p, q, \neg p\To p$}

\AxiomC{}
\RightLabel{$Id$}
\UnaryInfC{$q\To q$}
\RightLabel{$LW$}
\UnaryInfC{$p, q, \neg p\To q$}
\RightLabel{$R\land$}
\BinaryInfC{$p,q, \neg p\To p\land q$}

\AxiomC{}
\RightLabel{$Id$}
\UnaryInfC{$\neg p\To \neg p$}
\RightLabel{$LW$}
\UnaryInfC{$p,q, \neg p\To \neg p$}
\RightLabel{$R\neg\land$}
\UnaryInfC{$p,q, \neg p\To \neg (p\land q)$}

\RightLabel{$EXP_1^\con$}
\BinaryInfC{$\con(p\land q), p,q, \neg p\To$}
\end{prooftree}
\end{rmrk}

\section{Cut--elimination and applications} \label{sec.cut.elim}

In this section we shall prove the admissibility of the cut rule in 
the system \gletj.  Recall that a \textit{literal} is an atomic formula  or a negated atomic formula.

\begin{df}  \ 
\begin{enumerate} 
\item
The \textit{weight} of a formula is defined   as follows: 

	\begin{enumerate}%
		\item  $ w(l) = 0$, where \textit{l} is a literal;
            \item $w(\neg \neg A) =  w(A)  + 1$;
		\item  $w(A \ast B) =   w(A) +  w(B)  + 1$, for $\ast \in \{ \to,\land, \lor \}$;
		\item $w(\neg(A \ast B)) =   w(\neg A) +  w(\neg B)  + 1$,  
	 for $\ast \in \{ \to,\land, \lor \}$;
\item $w(\con A) =   w(A) +  w(\neg A)  + 1$;  
\item $w(\neg{\con}A) =   w(\con A)  + 1$. 

		\end{enumerate}	 
		




\item \label{gen.subf} {\it Generalized  subformula} 

\begin{enumerate} %
\item $A$ is 
a generalized subformula of $\neg A$;
\item $A$ and $B$ are    generalized subformulas of both  $A \ast B$,     
$\ast \in \{ \to, \lor, \land \}$;
\item $\neg A$ and $\neg B$ are    generalized subformulas of  $\neg (A \ast B)$, 
$\ast \in \{ \to, \lor, \land \}$;  
\item $A$ and $\neg A$ are    generalized subformulas of $\cons A $. 

\end{enumerate}
\end{enumerate}
\end{df}







In the sequel, we shall prove that \gletj\ enjoys the cut-elimination property, 
following the proof given by Gentzen (see \cite{Ta}). 
It is easy to check that the following version of the cut rule is equivalent to the one stated in \gletj .

\begin{prooftree}
\AxiomC{$\Gamma \To A$}
\AxiomC{$\Pi \To C$}
\RightLabel{\small $cut$, $A$}
\BinaryInfC{$\Gamma, \Pi_{A}  \To C $}
\end{prooftree}

\noindent where $A \in \Pi$ and $\Pi_{A}=\Pi - \{ A\}$; $A$ is said to be the {\em cut formula}.

Let $\P$ be a proof which contains only a cut in the last inference, we refer the left and right upper sequents as $S_1$ and $S_2$, respectively, and to the lower sequent as $S$. 

\begin{prooftree}
\AxiomC{$\vdots$}
\UnaryInfC{$S_1$}
\AxiomC{$\vdots$}
\UnaryInfC{$S_2$}
\RightLabel{\small $cut$, $A$}
\BinaryInfC{$S$}
\end{prooftree}

\mh We call a {\em left thread} of $\P$ ({\em right thread}) to any thread of $\P$ which contains the sequent $S_1$ ($S_2$). As usual, the {\em rank} of a left (right) thread is the number of consecutive sequents (counting upward from $S_1$ ($S_2$)) which contains the mix formula in its succedent (antecedent). We call $rank_l(\P)$ ($rank_r(\P)$) to the maximum of all ranks of left (right) threads in $\P$; and let $rank(\P)=rank_l(\P)+rank_r(\P)$. Note that $rank(\P)\geq 2$. Besides, we call the weight of $\P$, noted $w(\P)$, to the weight of the mix formula.  \\[2mm]

\begin{lm} If $\P$ is a proof  of the sequent $S$ which contains only a cut in the last inference, then there is a proof of $S$ without any cut.
\end{lm}
\begin{proof} We use double induction on the weight $n$, and the rank $m$ of $\P$ (i.e., transfinite induction on $\omega \cdot n + m$).\\[2mm]
{\bf Case I }: $m=2$. That is to say, $rank_{l}(\P)= 1 = rank_{r}(\P)$.\\[1.5mm]
Let $S_1$ and $S_2$ the left and right upper sequent of the mix in $\P$. We shall analyze different subcases. \\
{\bf (a)} If $S_1$ (or $S_2$) is an initial sequent (an instance of ($Id$)), it is clear that the theorem holds. 

\begin{prooftree}
\AxiomC{$A \To A$}

\alwaysNoLine
\AxiomC{$\vdots$}
\UnaryInfC{$\Pi \To C$}
\alwaysSingleLine
\RightLabel{\small $cut, A$}
\BinaryInfC{$A, \Pi_{A}  \To C$}
\end{prooftree}

\noindent This proof can be replaced by the next cut-free proof ($\alpha, \Pi_{A}$ is $\Pi$, since $A\in \Pi$)

\begin{prooftree}
\AxiomC{$\vdots$}
\alwaysNoLine
\UnaryInfC{$\Pi \To C$}
\end{prooftree}

\noindent {\bf (b)} $S_2$ is an initial sequent. Analogous to {\bf (a)}.\\
{\bf (c)} Then, suppose that neither $S_1$ nor $S_2$ is an initial sequent. If $S_1$ (or $S_2$) is the lower sequent of weakening, it is easy to check that the theorem holds.
Since $rank_{l}(\P)= 1$, $\alpha$ cannot appear in the succedent of the upper sequent of the weakening, therefore,  $\P$ has the form

\begin{prooftree}
\alwaysNoLine
\AxiomC{$\vdots$}
\UnaryInfC{$\Gamma \To $}
\alwaysSingleLine
\RightLabel{\small $RW$}
\UnaryInfC{$\Gamma \To  A$}

\alwaysNoLine
\AxiomC{$\vdots$}
\UnaryInfC{$\Pi \To C$}
\alwaysSingleLine
\RightLabel{\small $cut$}
\BinaryInfC{$\Gamma, \Pi_{A}  \To C$}

\end{prooftree}

\noindent Then, we can eliminate this cut as follows: 

\begin{prooftree}
\AxiomC{$\vdots$}
\noLine
\UnaryInfC{$\Gamma \To $}
\RightLabel{\small$LW$'s - $RW$ }
\UnaryInfC{$\Gamma, \Pi_{A}  \To C$}
\end{prooftree}

\noindent {\bf (d)} Suppose that both $S_1$ and $S_2$ are the lower sequents of logical inferences. Since $rank_{l}(\P)= 1 = rank_{r}(\P)$ we know that the mix formula on each side must be the principal formula of the logical inference. Let $A$ be the mix formula of $\P$. If $A$ is of the form  $B\wedge D$ or $B \vee D$ or $B\to D$ the proof goes as in the classical case. \\
Suppose that $A$ is of the form $\neg(B\lor D)$, then $\P$ is as follows.
\begin{prooftree}
\AxiomC{$\vdots$}
\UnaryInfC{$\Gamma \To \neg B$}
\AxiomC{$\vdots$}
\UnaryInfC{$\Gamma \To \neg D$}
\RightLabel{$R\neg\vee$}
\BinaryInfC{$\Gamma \Rightarrow \neg(B\lor D)$}

\AxiomC{$\vdots$}
\UnaryInfC{$\neg B, \neg D, \Pi\Rightarrow C$}
\RightLabel{$L\neg\vee$}
\UnaryInfC{$\neg (B\vee D), \Pi\Rightarrow C$}

\RightLabel{\small $cut$, $\neg (B\vee D)$}
\BinaryInfC{$\Gamma, \Pi \Rightarrow C$}
\end{prooftree}

\noindent where $\neg (B\vee D) \not\in \Pi$. Then, we can construct the following proof $\P_1$:

\begin{prooftree}
\AxiomC{$\vdots$}
\UnaryInfC{$\Gamma \To \neg B$}
\AxiomC{$\vdots$}
\UnaryInfC{$\neg B, \neg D, \Pi\To C$}
\RightLabel{\small $cut$, $\neg B$}
\BinaryInfC{$\Gamma, \neg D, \Pi_{\neg B} \To C$}
\end{prooftree}

\noindent Since $w(\P_1)<w(\P)$, by (I.H.) we know that there exists a cut-free proof $\P_1'$ of $\Gamma, \neg D, \Pi_{\neg B} \Rightarrow C$. Then, in turn, we can construct the following proof $\P_2$:
 
\begin{prooftree}
\AxiomC{$\vdots$}
\UnaryInfC{$\Gamma \To \neg D$}

\AxiomC{$\P_1'$}
\noLine
\UnaryInfC{$\vdots$}
\UnaryInfC{$\Gamma, \neg D, \Pi_{\neg B} \To C$}
\RightLabel{\small $cut$, $\neg D$}
\BinaryInfC{$\Gamma, \Gamma_{\neg D}, (\Pi_{\neg B})_{\neg D} \To C$}
\end{prooftree}

\noindent That is, $\P_2$ is a proof of the sequent \, $\Gamma, \Pi - \{\neg B, \neg D\} \To C$. Since  $w(\P_2)<w(\P)$, by (I.H), there exists a cut-free proof $\P_3$ of \, $\Gamma, \Pi - \{\neg B, \neg D\} \To C$, and by means of weakening, there exists a cut-free proof of \, $\Gamma, \Pi \To C$.\\[2mm]
If $A$ is of the form $\neg(B\land D)$, $\neg(B\to D)$ or $\neg\neg B$, the proof is analogous. Note that if $A$ has the form $\con B$, then the sequent $S_1$ has to be the lower sequent of $RW$, since no logical rule in \gletj\ introduce $\con$ in the succedent; and this case has been already analyzed in {\bf (c)}. This same happens if $A$ is of the form $\neg\con B$.\\[2mm]
{\bf Case II }: $m>2$. Then, $rank_{l}(\P)> 1$ \, and /or \, $rank_{r}(\P)>1$. Here, $S_1$ is $\Gamma\Rightarrow\Delta$ and $S_2$ is $\Pi\Rightarrow\Lambda$. We analyze different subcases. \\[1.5mm]
(I.1) $rank_{r}(\P)> 1$ 

 (I.1.1) \,  If  $A \in \Gamma$ or $A$ is $C$, then the construction of the desired proof is immediate. Indeed, we construct the desired proof as follows, respectively.

\begin{center}
\AxiomC{if $A\in \Gamma$}
\noLine
\UnaryInfC{$\vdots$}
\noLine
\UnaryInfC{$\Pi\To C$}
\UnaryInfC{$A, \Pi_{A} \To C$}
\RightLabel{\small $LW$'s}
\UnaryInfC{$\Gamma, \Pi_{A}\To C$}
\DisplayProof \hspace{1cm}
\AxiomC{if $A$ is $C$}
\noLine
\UnaryInfC{$\vdots$}
\noLine
\UnaryInfC{$\Gamma\To A$}
\UnaryInfC{$\Gamma\To C$}
\RightLabel{\small $LW$'s}
\UnaryInfC{$\Gamma, \Pi_{A}\To C$}
\DisplayProof
\end{center}

(I.1.2) \,  If $S_2$ is the lower sequent of an application of a rule $(r)$ which is not a logical rule with main formula $A$.  Then, $\P$ looks like:
\begin{prooftree}
\AxiomC{$\vdots$}
\UnaryInfC{$\Gamma \Rightarrow A$}

\AxiomC{$\vdots$}
\UnaryInfC{$\Pi' \Rightarrow C'$}
\RightLabel{\small $(r)$}
\UnaryInfC{$\Pi\Rightarrow C$}

\RightLabel{\small $cut$, $A$}
\BinaryInfC{$\Gamma, \Pi_{A} \Rightarrow C$}
\end{prooftree}
\noindent Since $rank_{r}(\P)> 1$, $A \in \Pi'$. Then, the following proof $\P$
\begin{prooftree}
\AxiomC{$\vdots$}
\UnaryInfC{$\Gamma \Rightarrow A$}

\AxiomC{$\vdots$}
\UnaryInfC{$\Pi' \Rightarrow C'$}
\RightLabel{\small $cut$, $A$}
\BinaryInfC{$\Gamma, \Pi'_{A} \Rightarrow C'$}
\end{prooftree}
contains only a cut in the last inference and $w(\P)=w(\P')$ but $rank(\P)<rank(\P')$. (Indeed, $rank_d(\P')=rank_d(\P)-1$). By the (I.H.), there is a cut-free proof of the sequent $\Gamma,\Pi'_{A}\To C'$. Then, using the rule $(r)$  we can construct the desired proof. 

(I.1.3) \, We can assume that $A \notin \Gamma$ and $S_2$ is the lower sequent of a logical inference whose principal formula is $A$.  If the $A$ is of the form $B \land D$,  $B\lor D$ or $B \to D$ the proof goes as in the classical case. Suppose that $A$ is of the form $\neg (B \to D)$, then $\P$ looks like

\begin{prooftree}
\AxiomC{$\vdots$}
\UnaryInfC{$\Gamma \To \neg(B\to D)$}

\AxiomC{$\vdots$}
\UnaryInfC{$\Pi, B, \neg D \To C$}
\RightLabel{$L\neg\to$}
\UnaryInfC{$ \Pi, \neg(B\to D) \To C$}

\RightLabel{\small $cut$, $ \neg(B\to D)$}
\BinaryInfC{$\Gamma, \Pi_{ \neg(B\to D)} \To C$}
\end{prooftree}

\noindent Since $rank_{r}(\P)> 1$,  $ \neg(B\to D) \in \Pi$. Then, we can construct proof $\P_1$.

\begin{prooftree}
\AxiomC{$\vdots$}
\UnaryInfC{$\Gamma \To \neg(B\to D)$}
\AxiomC{$\vdots$}
\UnaryInfC{$\Pi, B, \neg D \To C$}
\RightLabel{\small $cut$, $\neg(B\to D)$ }
\BinaryInfC{$\Gamma, \Pi_{ \neg(B\to D)}, B, \neg D \To C$}
\end{prooftree}

\noindent Clearly, $w(\P_1)=w(\P)$ but $rank(\P_1)=rank(\P)-1$. By the (I.H.), there is a cut-free proof $\P'_1$ of the sequents $\Gamma, \Pi_{ \neg(B\to D)}, B, \neg D \To C$. Then, consider the following proof $\P_2$

\begin{prooftree}
\AxiomC{$\vdots$}
\UnaryInfC{$\Gamma \To \neg(B\to D)$}

\AxiomC{$\P'_1$}
\noLine
\UnaryInfC{$\vdots$}
\UnaryInfC{$\Gamma, \Pi_{ \neg(B\to D)}, B, \neg D \To C$}
\RightLabel{$L\neg\to$ }
\UnaryInfC{$\Gamma, \Pi_{ \neg(B\to D)}, \neg(B\to D) \To C$}

\RightLabel{\small $cut$, $\neg(B\to D)$}
\BinaryInfC{$\Gamma, \Gamma, \Pi_{ \neg(A\to B)} \To C$}
\end{prooftree}

\noindent Then, $w(\P_2)=w(\P)$ but $rank_d(\P'')=1<rank_{d}(\P)$. By the (I.H.), there is a cut-free proof of $\Gamma, \Pi_{ \neg(A\to B)} \To C$.\\[2mm]
(I.2) $rank_{r}(\P)= 1$ and the $rank_{l}(\P)> 1$. This case is proved in the same way as (I.1) above.
\end{proof}

\begin{tr}\label{ElimCorte} \,  \, 

\m\mh \gletj\ admits cut-elimination.
\end{tr}
\begin{proof} By induction on the number of applications of the cut rule in a given \letj -proof.
\end{proof}

\noindent Some consequences of the cut-elimination theorem are the following.

\begin{cor}(Generalized subformula property) 

\m\mh In a cut-free proof $\P$ in \gletj\  all the formulas which occur in it are generalized subformulas of formulas occurring in the end sequent.  
\end{cor} 
\begin{proof}
By induction on the number of inferences in the cut-free proof $\P$.
\end{proof}

\begin{cor}\label{decidable} 
\m\mh \letj\  is decidable.
\end{cor}

\section{Sequent calculus for \letf} \label{sec.sequent.letf}

In this section we shall introduce a sequent calculus for \letf\ in the language $\{ \vee, \wedge, \neg, \con \}$. 

\begin{df} \ 

\m\mh The sequent calculus \gletf\ consists of the following axiom and rules.

 \

\textit{Axioms:}
\begin{center}
$\infer[Id] {A\To A} {}$

\end{center}

\textit{Structural rules:}
\begin{center}
$\infer[LW] {A, \Gamma \To \Delta} {\Gamma \To \Delta}$
\quad 
$\infer[RW] {\Gamma \To A, \Delta } { \Gamma\To\Delta}$

\

$\infer[Cut] {\Gamma \To \Delta} { \Gamma \To A, \Delta & A, \Gamma \To \Delta}$
\end{center}

\textit{Logical rules:}
\begin{center}
$\infer[L\land] {A \land B,\Gamma \To \Delta} {A,B , \Gamma \To \Delta}$  
\quad 
$\infer[R\land] {\Gamma \To A \land B, \Delta} { \Gamma\To A, \Delta   &   \Gamma\To B, \Delta}$

\

$\infer[L\lor] {A \lor B,\Gamma \To \Delta} {A,\Gamma \To \Delta & B,\Gamma \To \Delta}$
\quad 
$\infer[R\lor] {\Gamma \To A \lor B, \Delta } { \Gamma\To A, B, \Delta }$

\

$\infer[L\neg\lor] {\neg(A \lor B),\Gamma \To \Delta} {\neg A, \neg B,\Gamma \To \Delta}$
\quad 
$\infer[R\neg\lor] {\Gamma \To \neg(A \lor B), \Delta} { \Gamma\To \neg A, \Delta  &   \Gamma\To \neg B, \Delta}$

\

$\infer[L\neg\land] {\neg(A \land B),\Gamma \To \Delta} 
{\neg A,\Gamma \To C & \neg B,\Gamma \To \Delta}$
\quad 
$\infer[R\neg\land] {\Gamma \To \neg(A \land B), \Delta} {\Gamma\To \neg A, \neg B, \Delta}$

\

$\infer[L\con_1] {\con A, \Gamma\To \Delta}{\Gamma\To A, \Delta & \Gamma\To \neg A, \Delta}$
\quad 
$\infer[L\con_2] {\cons A, \Gamma \To \Delta}{A, \Gamma\To \Delta & \neg A, \Gamma\To \Delta }$

\

$\infer[L\neg\neg] {\neg\neg A , \Gamma \To  \Delta}{A , \Gamma \To \Delta}$
\quad 
$\infer[R\neg\neg] {\Gamma \To \neg\neg A, \Delta} {\Gamma \To A, \Delta}$

\

$\infer[L\neg\con]  { \neg\con A, \Gamma \To\Delta} {\Gamma \To A, \Delta}$
\quad 
$\infer[R\neg\con] { \Gamma \To \neg\con A , \Delta}{A , \Gamma \To \Delta}$

\end{center}
\end{df}

Note that \gletf\ has left and right rules for the negation connective $\neg$ interacting with each connective and  that the weight of the lower sequent of every rule us less that the weight of the upper(s) sequent.  This is a desired property for a syntactic proof of the cut elimination property. So,  the cut-elimination  
property for \gletf\ can be proved in a similar way to what was done for \letj. 

\begin{tr} \,  \, 

\m\mh \gletf\ admits cut-elimination.
\end{tr}
\begin{proof}
Analogous to the case \gletj. 
\end{proof}

\begin{cor} (Generalized subformula property) 

\m\mh In a cut-free proof $\P$ in \gletf\  all the formulas which occur in it are generalized subformulas of formulas occurring in the end sequent.  
\end{cor}

\begin{cor}\label{decidableLetj} \,  \, 

\m\mh \letf\  is decidable.
\end{cor}

\section{Final remarks}

The logics \letj\ and \letf\ were originally presented in \cite{letj} and \cite{letf} 
in the form of natural deduction systems. 
The main reason for this was to emphasize that the intended approach was syntactic rather than semantic, 
since the motivation was to formalize argumentative contexts in which people make inferences with evidence, 
both conclusive and non-conclusive \cite[pp.~3793f.]{letj}. 
In \cite[Sect.~2.2.1]{letf} this intuitive interpretation was extended to reliable and unreliable information. 
 In addition, the semantics proposed so far for \letj\ and \letf\ are non-compositional, 
 and therefore unable to provide adequate explanations of the meanings of their formulas.    
 The prospect of explaining the meanings of the expressions of \letj\ by means of an inferential semantics 
 was indeed mentioned in \citep{Carnielli2019a}, which makes obtaining normalization all the more desirable. 
 
 Presenting  natural deduction systems without the corresponding normalization theorems 
 has indeed the flavor of an unfinished task.  
 In this paper, the first strictly dedicated to the proof theory of \lets, we 
 fill this gap and also pave the way for further investigations, 
 including the proof theory of other logics equipped with recovery operators. 

Decision procedures for the logic \letf\ were already provided by its valuation semantics \cite{letf} and by analytic tableaux \cite{tableaux.lets},  but a decision procedure for \letj\ had not yet been presented. Decidability by means of sequent calculus  also fills this gap. 

Normalization theorems for the \con-free fragments of \letj\ and \letf, 
namely \nel\ and \fde, have already been presented (see e.g. \citep{Kurbis_Petrukhin_2021}), 
as well as cut-elimination theorems (see e.g. \citep{kamide.wansing.2015}).  
 The central point here   was how to deal with the rules for the operator \con. 
 More precisely, the treatment given by Kamide in~\cite{Kam} for explosion and excluded middle rules 
 inspired us to obtain a suitable normalization theorem. 
 On the other hand, its is worth noting that a direct translation of the natural deduction version 
 to sequent system of both logics produced systems which do not enjoy the cut-elimination property. 
  The corresponding rules for explosion and excluded middle laws in a cut-free sequent system turn out to be 
 left introduction rules of the \con\ connective. 
  Note that this fits with the fact that there is no introduction rule for \con\ in the natural deduction systems, 
 which is motivated by the ideia that the information about realiability of a formula comes always from outside the formal system.

\bibliographystyle{plain}

%
%
%
%
%

\end{document}